\documentclass[11pt]{amsart}

\usepackage{mathtools,amsmath,amsthm,amsfonts,latexsym,amssymb,mathrsfs,setspace,verbatim,mathtools,comment,enumitem,doi}
\usepackage[square,sort,comma,numbers]{natbib}
\usepackage{fullpage}
\usepackage{hyperref}

\mathtoolsset{showonlyrefs=true}

\newtheorem{thm}{Theorem}

\newtheorem{lem}[thm]{Lemma}

\newtheorem{cor}[thm]{Corollary}
\numberwithin{thm}{section}
\numberwithin{equation}{section}

\theoremstyle{definition}

\newcommand{\rat}{\mathbb Q}
\newcommand{\real}{\mathbb R}

\newcommand{\alg}{\overline\rat}
\newcommand{\algt}{\alg^{\times}}

\newcommand{\intg}{\mathbb Z}
\newcommand{\nat}{\mathbb N}

\newcommand{\tors}{\mathrm{tors}}
\newcommand{\supp}{\mathrm{supp}}
\newcommand{\spann}{\mathrm{span}}

\newcommand{\R}{\mathcal R}
\newcommand{\A}{\mathcal A}

\newcommand{\generalring}{\mathfrak D}

\begin{document}

\title[Consistent maps and measures]{The correspondence between consistent maps and measures on the places of $\overline{\mathbb Q}$}

\author[C.L. Samuels]{Charles L. Samuels}
\address{Christopher Newport University, Department of Mathematics, 1 Avenue of the Arts, Newport News, VA 23606, USA}
\email{charles.samuels@cnu.edu}
\thanks{This work was funded in part by the AMS-Simons Research Enhancement Grant for PUI Faculty}

\begin{abstract}
	Recent work of the author established dual representation theorems for certain vector spaces that arise in an important article of Allcock and Vaaler.  These results constructed an object called
	a consistent map which acts like a measure on the set of places of $\overline{\mathbb Q}$, but is not a Borel measure on this space.  We describe the appropriate ring of sets $\mathcal R$ for which every consistent 
	map arises from a measure on $\mathcal R$.  We further obtain the conditions under which a consistent map may be extended to a measure on the smallest algebra containing $\mathcal R$.	
\end{abstract}

\maketitle

\section{Introduction}\label{Intro}

Let $\alg$ be a fixed algebraic closure of $\rat$ and let $\algt_\tors$ denote the group of roots of unity.  We observe that $\mathcal G := \algt/\algt_\tors$ is a vector space over $\rat$
with addition and scalar multiplication defined by
\begin{equation*}
	(\alpha,\beta)\mapsto \alpha\beta\quad\mbox{and}\quad (r,\alpha) \mapsto \alpha^r.
\end{equation*}
For each number field $K$, we write $M_K$ to denote the set of all places of $K$.   If $L/K$ is a finite extension and $w\in M_L$, then $w$ divides a unique place $v$ of $K$, and in this case, 
we write $K_w$ to denote the completion of $K$ with respect to $v$.  Additionally, if $v$ divides the place $p$ of $\rat$, then we let $\|\cdot \|_v$ be the unique extension to $K_v$ of the usual $p$-adic absolute 
value on $\rat_v$

Let $Y$ denote the set of all places of $\alg$ and define $Y(K,v) = \{y\in Y:y\mid v\}$.  Further setting
\begin{equation} \label{Jdef}
	\mathcal J = \left\{(K,v):[K:\rat]<\infty,\ v\in M_K\right\}\quad\mbox{and}\quad Y(\mathcal J) = \{Y(K,v): (K,v)\in \mathcal J\},
\end{equation}
Allcock and Vaaler \cite{AllcockVaaler} observed that $Y(\mathcal J)$ is a basis for a totally disconnected, Hausdorff topology on $Y$, and moreover, there is a Borel measure $\lambda$ on $Y$ such that
\begin{equation} \label{LambdaMeasure}
	\lambda(Y(K,v)) = \frac{[K_v:\rat_v]}{[K:\rat]}
\end{equation}
for all $(K,v)\in \mathcal J$.
Each element $\alpha\in \mathcal G$ corresponds to a locally constant function $f_\alpha:Y\to \real$ given by the formula $$f_\alpha(y) = \log \|\alpha\|_y.$$  When $\mathcal G$ is equipped with a
norm arising from the Weil height, Allcock and Vaaler \cite{AllcockVaaler} proved that $\alpha\mapsto f_\alpha$ defines an isometric isomorphism from $\mathcal G$ onto a dense $\rat$-linear subspace of
\begin{equation*}
	\left\{f\in L^1(Y): \int_Y f(y)d\lambda(y) = 0\right\}.
\end{equation*}

More recently, the author \cite{SamuelsClassification,SamuelsConsistent,SamuelsRep} began the study of various function spaces on $\mathcal G$ such as the collection $\mathcal L(\mathcal G,\real)$ of $\rat$-linear maps from
$\mathcal G$ to $\real$.  Spaces of this sort are important in number theory because they often classify extensions of completely additive arithmetic functions to linear maps on $\mathcal G$.  
In attempting to describe these spaces, the Riesz Representation Theorem \cite[Theorems 2.14 \& 6.19]{Rudin} is a standard starting point, establishing a correspondence between continuous linear maps $\mathcal G\to \real$ and 
signed Borel measures on $Y$.  Nevertheless, if we seek a classification of all (not necessarily continuous) linear maps from $\mathcal G$ to $\real$ or all continuous linear maps from $\mathcal G$ to $\rat$, this famous result
fails to address these goals.  

In \cite{SamuelsClassification,SamuelsConsistent,SamuelsRep}, we resolved this discrepancy by considering the following object.
A map $c:\mathcal J\to \real$ is called {\it consistent} if 
\begin{equation*}
	c(K,v) = \sum_{w\mid v} c(L,w)
\end{equation*}
for all $(K,v)\in \mathcal J$ and all finite extensions $L/K$.  We shall write $\mathcal J^*$ to denote the vector space of all consistent maps, where addition and scalar multiplication are given by the formulas
\begin{equation} \label{ConsistentOperations}
	(c+d)(K,v) = c(K,v) + d(K,v) \quad\mbox{and}\quad (rc)(K,v) = rc(K,v).
\end{equation}
One basic example is given by the right hand side of \eqref{LambdaMeasure} which we shall simply denote by $\lambda$.  For each $c\in \mathcal J^*$ and $\alpha\in \alg^\times$, we let
\begin{equation} \label{PhiDef}
	\Phi_c(\alpha) = \sum_{v\in M_K} c(K,v)\log\|\alpha\|_v,
\end{equation}
where $K$ is any number field containing $\alpha$.  Because $c$ is assumed to be consistent, this definition does not depend on the choice of $K$, and moreover, its value is unchanged when $\alpha$ is multiplied by a root of unity.

The main result of \cite{SamuelsRep} shows that $c\mapsto \Phi_c$ is a surjective $\real$-linear transformation from $\mathcal J^*$ to $\mathcal L(\mathcal G,\real)$ whose kernel is equal to $\spann_\real\{\lambda\}$.
If we further write $$\mathcal J_\infty^* = \{c\in \mathcal J^*: c(\rat,\infty) = 0\},$$ then $c\mapsto \Phi_c$ is an $\real$-vector space isomorphism from $\mathcal J_\infty^*$ onto $\mathcal L(\mathcal G,\real)$.
As is noted in \cite{SamuelsRep}, these results provide a vehicle to describe extensions of completely additive functions $\nat\to \real$ to elements of $\mathcal L(\mathcal G,\real)$.  For instance, the {\it prime Omega function}
is defined so that $\Omega(n)$ is the number of (not necessarily distinct) prime factors of $n$.  Assuming that $c\in \mathcal J_\infty^*$, 
then $\Phi_c:\mathcal G\to \real$ is an extension of $\Omega$ if and only if $c(\rat,p) = -1/\log p$ for all primes $p\ne \infty$.

As the reader may have observed, the summation on the right hand side of \eqref{PhiDef} heavily mirrors the definition of an integral, where $c$ plays the role of a measure.  In fact, if $\mu$ is a signed Borel measure 
which is finite on compact subsets of $Y$ then
\begin{equation} \label{MeasureToConsistent}
	c(K,v) := \mu(Y(K,v))
\end{equation}
defines a consistent map, and furthermore, the right hand side of \eqref{PhiDef} is precisely equal to 
\begin{equation*}
	\int_{Y} f_\alpha(y)d\mu(y).
\end{equation*}
We used similar methods in \cite{SamuelsClassification,SamuelsConsistent} to obtain isomorphisms with the duals of $\algt/\overline{\intg}^\times$ and $LC_c(Y)$, respectively.  While the precise
summations are slightly different than \eqref{PhiDef}, they again borrow heavily from principles of integration theory.  

All of the above observations suggest that consistent maps really are a sort of measure and the right hand side of \eqref{PhiDef} may always be expressed as an integral.
While these are tantalizing possibilities, signed Borel measures on $Y$ are not the correct objects to consider.  Indeed, the work of Aberg and the author \cite{AbergSamuels} created a counterexample showing 
that not all consistent maps arise from signed Borel measures in the manner given by \eqref{MeasureToConsistent}.  

We let $\R$ be the smallest ring of sets on $Y$ containing $Y(\mathcal J)$, or in other words, $\R$ is the collection of all finite unions of sets in $Y(\mathcal J)$ along with the empty set.
The following is the first of our two main results.

\begin{thm} \label{RingMeasure}
	If $c:\mathcal J\to \real$ is a consistent map then there exists a unique signed measure $\mu$ on $\R$ such that $c(K,v) = \mu(Y(K,v))$ for all $(K,v)\in \mathcal J$.  In this case, $f_\alpha$ is integrable with
	respect to $\mu$ and 
	\begin{equation*}
		\Phi_c(\alpha) = \int_Yf_\alpha(y) d\mu(y)
	\end{equation*}
	for all $\alpha\in \mathcal G$.
\end{thm}

For each $\alpha\in \mathcal G$, we remark that $\supp(f_\alpha)\in \mathcal R$.  Moreover, if $K$ is a number field containing some coset representative of $\alpha$, then $f_\alpha$ is constant on $Y(K,v)$
for all $v\in M_K$.  These observations mean that the integral appearing in Theorem \ref{RingMeasure} is well-defined even though $Y\not\in \R$.

When combining this result with \cite{SamuelsRep}, we can express extensions of completely additive arithmetic functions as integrals.  For example, every extension of the prime Omega function
to a linear map $\Omega: \mathcal G\to \real$ is given by the formula
\begin{equation*}
	\Omega(\alpha) = \int_Yf_\alpha(y) d\mu(y),
\end{equation*}
where $\mu$ is a measure on $\R$ satisfying
\begin{equation*}
	\mu(Y(\rat,p)) = \begin{cases} -1/\log p & \mbox{if } v\mid p\mbox{ and } p \ne \infty \\
						0 & \mbox{if } v\mid \infty.\end{cases}
\end{equation*}
Other completely additive arithmetic functions can be handled in a similar manner, and in this way, we have learned to express such functions as integrals over $Y$.

For simplicity, we now write $c:\mathcal R\to \real$ to denote the measure described in Theorem \ref{RingMeasure}.  Since $Y$ cannot be expressed as a finite union of sets in $Y(\mathcal J)$, the collection 
$\mathcal R$ is certainly not an algebra of sets, and hence, we arrive at a natural follow-up question.  Which consistent maps have a further extension to the smallest algebra $\A$ on $Y$ containing $\mathcal R$?   

Importantly, not all consistent maps may be extended in such a way.  For example, let $p_1,p_2,p_3,\ldots$ be the complete list or rational primes (including $\infty$) and define
\begin{equation*}
	c(K,v) = (-1)^n \frac{[K_v:\rat_{p_n}]}{[K:\rat]}\quad\mbox{ for all } v\mid p_n.
\end{equation*}
Using \cite[Eq. (2)]{Frohlich}, it is straightforward to prove that $c$ is a consistent map.  On the other hand, we have 
\begin{equation*}
	Y = \bigcup_{n=1}^\infty Y(\rat,p_n),
\end{equation*}
while $\sum_{n=1}^\infty c(\rat,p_n)$ fails to converge unconditionally to a value in $[-\infty,\infty]$.  As a result, the required countable additivity property for measures cannot be satisfied.
In order to eliminate cases of this sort, we say that a subcollection $\Gamma \subseteq Y(\mathcal J)$ is a {\it basis partition of $Y$} if the elements of $\Gamma$ are pairwise disjoint and
\begin{equation*}
	Y = \bigcup_{A\in \Gamma} A.
\end{equation*}
A consistent map $c:\mathcal J\to \real$ is called {\it globally consistent} if $\sum_{A\in \Gamma} c(A)$ converges unconditionally in $[-\infty,\infty]$ for every basis partition $\Gamma$ of $Y$.
This property completely classifies those consistent maps that arise from measures on $\mathcal A$.

\begin{thm} \label{AlgebraMeasure}
	If $c:\mathcal J\to \real$ is a globally consistent map then there exists a unique signed measure $\mu$ on $\A$ such that $c(K,v) = \mu(Y(K,v))$ for all $(K,v)\in \mathcal J$.
\end{thm}

The remainder of this article is organized in the following way.  Our main results require some preliminary work on finitely additive measures, often called charges as in \cite{RaoRao}, which we present
in Section \ref{Rings}.  In Section \ref{Measures}, we provide a definition for a specific class of measures on $\R$, and we prove that they are in one-to-one correspondence with consistent maps.
Theorem \ref{RingMeasure} is obtained as a direct corollary to that result.  Finally, we use Section \ref{Algebra} to study globally consistent maps in greater detail and to establish Theorem \ref{AlgebraMeasure}.

\section{Rings and Charges} \label{Rings}

Following the definition presented in \cite[Ch. 1]{Kesavan}, a non-empty collection $\generalring$ of subsets of $Y$ is called a {\it ring of sets} if whenever $A,B\in \generalring$ then 
$A\cup B\in \generalring$ and $A\setminus B\in \generalring$.  Clearly $\emptyset\in \generalring$ for any ring $\generalring$, and moreover, if $A,B\in \generalring$ then $A\cap B\in \generalring$.  
If, in addition, $Y\in \generalring$, then $\generalring$ is called an {\it algebra of sets}.   We begin with the following basic facts about $\R$ as defined prior to the statement of Theorem \ref{RingMeasure}.

\begin{lem}\label{RFacts}
	The following properties hold:
	\begin{enumerate}[label=(\roman*)]
		\item\label{OpenCompact} Suppose $A\subseteq Y$.  Then $A\in \mathcal R$ if and only if $A$ is open and compact.
		\item\label{Closed} If $A\in \mathcal R$ then $A$ is closed.
		\item\label{Ring} $\R$ is the smallest ring of sets on $Y$ for which $Y(\mathcal J)\subseteq \R$.
		\item\label{PairwiseDisjoint} If $A$ is a non-empty element of $\R$ then there exists a pairwise disjoint collection $\{W_i:1\leq i\leq n\} \subseteq Y(\mathcal J)$ such that
			\begin{equation*}
				A = \bigcup_{i=1}^n W_i.
			\end{equation*}
		\item\label{InfinitePairwiseDisjoint} If $\{V_i: i\in \nat\} \subseteq Y(\mathcal J)$ is a pairwise disjoint collection with $V_i\ne \emptyset$ for all $i$ then 
			\begin{equation*}
				\bigcup_{i=1}^\infty V_i \not \in \R.
			\end{equation*}
	\end{enumerate}
\end{lem}
\begin{proof}
	{\bf \ref{OpenCompact}:} First assume that $A\in \R$.  If $A=\emptyset$ then clearly $A$ is both open and compact.  Otherwise, $A$ is a finite union of sets of the form $Y(K,v)$, so according to \cite{AllcockVaaler}, 
	it is a finite union of open and compact sets.
	Therefore, $A$ must also be open and compact.  Now suppose that $A$ is open and compact.  For each point $y\in A$, there exists a number field $K_y$ and a place $v_y$ of $K_y$ such that $y\in Y(K_y,v_y) \subseteq A$.
	This means that $\{Y(K_y,v_y):y\in A\}$ is an open cover of $A$.  Now the compactness of $A$ implies that $A$ is a finite union of sets of the form $Y(K,v)$, as is required to belong to $\R$.
	
	{\bf \ref{Closed}:} The work of \cite{AllcockVaaler} shows that $Y(K,v)$ is closed for all $(K,v)\in \mathcal J$.  Hence, \ref{Closed} follows from the fact that $A$ is a finite union of these sets.

	{\bf \ref{Ring}:} If $A,B\in \R$ then it follows directly from the definition of $\R$ that $A\cup B\in \R$.  By using \ref{OpenCompact} and \ref{Closed}, we have that $A\setminus B = A\cap (Y\setminus B)$ is a closed subset of a compact set, and hence, 
		is itself compact.  Of course, \ref{OpenCompact} and \ref{Closed} also imply that $A\setminus B$ is open so it must belong to $\R$.  These observations prove that $\R$ is indeed a ring of sets, and since each ring must be closed
		under finite unions, certainly $\R$ is the smallest ring containing $Y(\mathcal J)$.
		
	{\bf \ref{PairwiseDisjoint}:} We let $U_i\in Y(\mathcal J)$ be such that
		\begin{equation*}
			A = \bigcup_{i=1}^n U_i.
		\end{equation*}
		For each $i$, let $K_i$ be a number field and $v_i$ a place of $K_i$ such that $U_i = Y(K_i,v_i)$.  Letting $L$ be the compositum of the fields $K_i$, then each set $U_i$ is a disjoint union
		\begin{equation*}
			U_i = \bigcup_{w\mid v_i} Y(L,w),
		\end{equation*}
		and we obtain
		\begin{equation} \label{AlmostDisjoint}
			A = \bigcup_{i=1}^n \bigcup_{w\mid v_i} Y(L,w).
		\end{equation}
		If $w_1$ and $w_2$ are places of $K$, then either $Y(L,w_1) = Y(L,w_2)$ or $Y(L,w_1)\cap Y(L,w_2) = \emptyset$ so \ref{PairwiseDisjoint} now follows by removing duplicate terms from \eqref{AlmostDisjoint}.
	
	{\bf \ref{InfinitePairwiseDisjoint}:} Now let $B$ denote the countably infinite union described in \ref{InfinitePairwiseDisjoint}.  If $B\in \R$ then we know that $B$ is compact, but this contradicts the fact that $\{V_i: i\in \nat\}$ is an open cover
		of $B$ by a countably infinite disjoint collection of non-empty sets.
\end{proof}

Our proof of Theorem \ref{RingMeasure} requires that we first consider a finitely additive version of a measure, often called a charge.
Supposing that $\generalring$ is any ring of sets on $Y$, a map $\mu:\generalring\to \real$ is called a {\it charge on $\generalring$} if it satisfies the following conditions:
\begin{enumerate}[label=(C\arabic*)]
	\item\label{SimpleEmpty} $\mu(\emptyset) = 0$
	\item\label{FiniteAdditive} If $A,B\in \generalring$ are such that $A\cap B = \emptyset$ then $\mu(A\cup B) = \mu(A) + \mu(B)$.
\end{enumerate}
There is some existing study of charges in the literature (see \cite{RaoRao,YosidaHewitt,Kunisada}, for instance).  While their definitions of charge always enforce properties \ref{SimpleEmpty} and \ref{FiniteAdditive},
they typically differ from our definition in important ways. For example, \cite{YosidaHewitt} requires that charges be bounded while \cite{Kunisada} requires that they be both bounded and non-negative.  
As consistent maps need not satisfy either of these conditions, assumptions of this sort are incompatible with our goals.  Our next lemma assists with our task of extending a consistent map to a charge on $\mathcal R$.

\begin{lem} \label{ConsistentExtension}
	Let $c:\mathcal J\to \real$ be a consistent map and $m,n\in \nat$.  For each $1\leq i\leq m$ and $1\leq j\leq n$, let $V_i,W_j\in Y(\mathcal J)$ be such that $\{V_i:1\leq i\leq m\}$ and $\{W_j:1\leq i\leq n\}$
	are each pairwise disjoint.  If
	\begin{equation} \label{UnionEquality}
		\bigcup_{i=1}^m V_i = \bigcup_{j=1}^n W_j
	\end{equation}
	then
	\begin{equation*}
		\sum_{i=1}^m c(V_i) = \sum_{j=1}^n c(W_j)
	\end{equation*}
\end{lem}
\begin{proof}
	Let $(K_i,v_i),(L_j,w_j)\in \mathcal J$ be such that $V_i = Y(K_i,v_i)$ and $W_j = Y(L_j,w_j)$ and let $E$ be the compositum of the fields $K_i$ and $L_j$.  For each place $u\in M_E$, we claim that the following are equivalent:
	\begin{enumerate}[label={(\roman*)}]
		\item\label{v} $u\mid v_i$ for some $1\leq i\leq m$
		\item\label{vU} There exists a unique point $1\leq i\leq m$ such that $u\mid v_i$
		\item\label{w} $u\mid w_j$ for some $1\leq j\leq n$
		\item\label{wU} There exists a unique point $1\leq j\leq n$ such that $u\mid w_j$.
	\end{enumerate}
	To see that \ref{v}$\implies$\ref{vU}, we assume that $u$ divides $v_j$ for some $1\leq j\leq m$.  This means that the non-empty set $Y(E,u)$ is a subset of both $Y(K_i,v_i)$ and $Y(K_j,v_j)$ contradicting our assumption that $V_i$ are pairwise disjoint.
	A similar argument shows that \ref{w}$\implies$\ref{wU}.  For \ref{vU}$\implies$\ref{w}, our assumption means that $Y(E,u)\subseteq Y(K_i,v_i) = V_i$.  Then using \eqref{UnionEquality}, we conclude that 
	\begin{equation*}
		Y(E,u)\subseteq \bigcup_{j=1}^n W_j = \bigcup_{j=1}^n Y(L_j,w_j).
	\end{equation*}
	Since $Y(E,u)$ is non-empty, there must exist $j$ such that $Y(E,u)\cap Y(L_j,w_j) \ne\emptyset$, and it follows that $u\mid w_j$.  A similar argument shows that \ref{wU}$\implies$\ref{v}, finally showing that the above four statements
	are equivalent.  As a result, we now conclude that
	\begin{equation*}
		\sum_{i=1}^m \sum_{u\mid v_i} c(E,u) = \sum_{j=1}^n \sum_{u\mid w_j} c(E,u).
	\end{equation*}
	Using the consistency of $c$, we find that
	\begin{align*}
		\sum_{i=1}^m c(V_i) = \sum_{i=1}^m c(K_i,v_i) & = \sum_{i=1}^m \sum_{u\mid v_i} c(E,u) \\
			& = \sum_{j=1}^n \sum_{u\mid w_j} c(E,u) = \sum_{j=1}^n c(L_j,w_j) = \sum_{j=1}^n c(W_j),
	\end{align*}
	establishing the lemma.	
\end{proof}

One easily verifies that the collection $\mathcal M_{\generalring}^0$ of charges on $\generalring$ forms a vector space over $\real$ with addition and scalar multiplication given by
\begin{equation*}
	(\mu + \nu)(A) = \mu(A) + \nu(A)\quad\mbox{and}\quad (r\mu)(A) = r\mu(A).
\end{equation*}
If $\mu\in \mathcal M_\R^0$ and $(K,v)\in \mathcal J$, we shall write 
\begin{equation*}
	c_\mu(K,v) = \mu(Y(K,v)).
\end{equation*}
Lemma \ref{ConsistentExtension} allows us to use this map to prove that $\mathcal M_{\R}^0$ and $\mathcal J^*$ are isomorphic.

\begin{thm}\label{SimpleMeasureIsomorphism}
	The map $\mu\mapsto c_\mu$ is a vector space isomorphism from $\mathcal M_{\R}^0$ to $\mathcal J^*$.
\end{thm}
\begin{proof}
	For each $\mu\in \mathcal M_{\R}^0$, we must first show that $c_\mu$ is a consistent map.  To see this, assume that $K$ is a number field, $L/K$ is a finite extension, and $v\in M_K$.  Then we have the disjoint union
	\begin{equation*}
		Y(K,v) = \bigcup_{w\mid v} Y(L,w)
	\end{equation*}
	so property \ref{FiniteAdditive} implies that
	\begin{equation*}
		\sum_{w\mid v} c_{\mu}(L,w) = \sum_{w\mid v} \mu(Y(L,w)) = \mu(Y(K,v)) = c_\mu(K,v)
	\end{equation*}
	showing that $c_\mu$ is indeed a consistent map.  For the remainder of the proof we shall define $\phi:\mathcal M_{\R}^0\to \mathcal J^*$ by $\phi(\mu) = c_{\mu}$.
	
	Letting $\mu,\nu\in \mathcal M_{\R}^0$ and $(K,v)\in \mathcal J$, we have
	\begin{align*}
		\phi(\mu+\nu)(K,v) & = c_{\mu+\nu}(K,v) \\
			& = (\mu + \nu)(Y(K,v)) \\
			& = \mu(Y(K,v)) + \nu(Y(K,v)) \\
			& = \phi(\mu)(K,v) + \phi(\nu)(K,v) \\
			& = (\phi(\mu) + \phi(\nu))(K,v),
	\end{align*}
	implying that $\phi(\mu + \nu) = \phi(\mu) + \phi(\nu)$.  Using a similar argument, we obtain that $\phi(r\mu) = r\phi(\mu)$, showing that $\phi$ is a linear transformation.
	
	To see that $\phi$ is injective, we let $\mu\in \mathcal M_{\R}^0$ and assume that $\phi(\mu) = 0$.  Now let $A\in \R$, so according to Lemma \ref{RFacts}\ref{PairwiseDisjoint},
	there exist number fields $K_i$ and places $v_i\in M_K$ such that $A$ is represented as the disjoint union
	\begin{equation*}
		A = \bigcup_{i=1}^n Y(K_i,v_i).
	\end{equation*}
	Therefore, 
	\begin{equation*}
		\mu(A) = \sum_{i=1}^n \mu(Y(K_i,v_i)) = \sum_{i=1}^n \phi(\mu)(K_i,v_i),
	\end{equation*}
	and since each summand on the right hand side is assumed to be $0$, we conclude that $\mu(A) = 0$, as required.
	
	To prove that $\phi$ is surjective, we assume that $c\in \mathcal J^*$ and must find $\mu\in \mathcal M_{\R}^0$ such that $c_\mu = c$.  As before, each non-empty set $A\in \R$ may be expressed as
	the disjoint union
	\begin{equation*}
		A = \bigcup_{i=1}^n Y(K_i,v_i),
	\end{equation*}
	for some number fields $K_i$ and places $v_i\in M_K$.  Now define $\mu:\R\to \real$ so that $\mu(\emptyset) = 0$ and
	\begin{equation} \label{muDef}
		\mu(A) = \sum_{i=1}^n c(K_i,v_i).
	\end{equation}
	In view of Lemma \ref{ConsistentExtension}, the right hand side of \eqref{muDef} is independent of the choice of number fields $K_i$ and places $v_i$.  Therefore, $\mu$ is a well-defined map
	from $\R$ to $\real$.  Clearly $c_\mu(K,v) = \mu(Y(K,v)) = c(K,v)$ so that $c_\mu = c$, and moreover, $\mu$ plainly satisfies property \ref{SimpleEmpty}.
		
	To prove \ref{FiniteAdditive}, let $A$ and $B$ be disjoint sets in $\R$.  Once again, we may express $A$ and $B$ as disjoint unions in the form
	\begin{equation*}
		A = \bigcup_{i=1}^n Y(K_i,v_i)\quad\mbox{and}\quad B = \bigcup_{j=1}^m Y(L_j,w_j).
	\end{equation*}
	Since $A$ and $B$ are themselves disjoint, we also have the disjoint union
	\begin{equation*}
		A\cup B = \left[\bigcup_{i=1}^n Y(K_i,v_i)\right] \bigcup \left[\bigcup_{j=1}^m Y(L_j,w_j)\right].
	\end{equation*}
	Therefore, \eqref{muDef} implies that
	\begin{equation*}
		\mu(A\cup B) = \sum_{i=1}^n c(K_i,v_i) + \sum_{j=1}^m c(L_j,w_j) = \mu(A) + \mu(B)
	\end{equation*}
	as required.
\end{proof}

\section{Measures} \label{Measures}

Although Theorem \ref{SimpleMeasureIsomorphism} shows how to represent consistent maps as charges, we still appear to be short of our stated goal proving that every consistent map arises from a measure on $\R$.
After all, a measure is required to satisfy the stronger countable additivity property rather than simply the finite additivity property \ref{FiniteAdditive}.  Luckily, the ring $\R$ is constructed in such a way that charges and certain 
types of measures are equivalent.  

Suppose $I$ is an arbitrary countable set and $\{x_i\}_{i\in I}$ is a collection of (not necessarily distinct) real numbers.  We say that
\begin{equation} \label{GeneralSeries}
	\sum_{i\in I} x_i
\end{equation}
{\it converges unconditionally} to a point $x\in [-\infty,\infty]$ if the series $\sum_{n=1}^\infty x_{\sigma(n)}$ converges to $x$ for all bijections $\sigma:\nat\to I$.
It is well-known that \eqref{GeneralSeries} converges unconditionally to real number if and only if it converges absolutely, and hence, the literature often uses these two concepts interchangeably.  
However, we must consider the possibility that $x = \pm \infty$ so there is a slight distinction between unconditional convergence and absolute convergence in our definition.

Our definition of unconditional convergence may be extended to sequences $\{x_i\}_{i\in I}$ having $x_i\in [-\infty,\infty]$.  Indeed, we may set
\begin{equation*}
	I_{0} = \{i\in I: |x_i|<\infty\},\quad I_{+} = \{i\in I:x_i = \infty\}\quad\mbox{and}\quad I_{-} = \{i\in I:x_i = -\infty\}
\end{equation*}
and consider three cases:
\begin{enumerate}[label=(U\arabic*)]
	\item\label{UCBasic} If $I_{+} = I_{-} = \emptyset$ and $\sum_{i\in I_0} x_i$ converges unconditionally to $x\in [-\infty,\infty]$, then we say that $\sum_{i\in I} x_i$ {\it converges unconditionally to $x$}.
	\item\label{UC+} If $I_{-} = \emptyset$, $I_{+} \ne\emptyset$ and $\sum_{i\in I_0} x_i$ converges unconditionally to $x\in (-\infty,\infty]$, then we say that $\sum_{i\in I} x_i$ {\it converges unconditionally to $\infty$}.
	\item\label{UC-} If $I_{+} = \emptyset$, $I_{-} \ne\emptyset$ and $\sum_{i\in I_0} x_i$ converges unconditionally to $x\in [-\infty,\infty)$, then we say that $\sum_{i\in I} x_i$ {\it converges unconditionally to $-\infty$}.
\end{enumerate}
If none of the above three conditions are satisfied, then $\sum_{i\in I} x_i$ fails to converge unconditionally, or simply {\it converges conditionally}.  
In the special case where $I$ is finite, unconditional convergence is not automatically satisfied. Under this assumption, $\sum_{i\in I} x_i$ converges unconditionally if and only if $I_+$ or $I_-$ is empty.   

Suppose that $\generalring$ is any ring of sets on $Y$.  We say that $\mu:\generalring\to [-\infty,\infty]$ is a {\it semi-finite signed measure on $\generalring$} if the following properties hold:
\begin{enumerate}[label={(M\arabic*)}]
	\item\label{ZeroEmpty} $\mu(\emptyset) = 0$
	\item\label{CompactFinite} $|\mu(A)| < \infty$ for all compact sets $A\in \generalring$
	\item\label{CountablyAdditiveList} If $\Gamma$ is a countable collection of pairwise disjoint sets in $\generalring$ such that $\bigcup_{A\in\Gamma} A\in \generalring$, then 
			$\sum_{A\in \Gamma} \mu(A)$ converges unconditionally to $\mu\left(\bigcup_{A\in \Gamma} A \right)$.
\end{enumerate}
We write $\mathcal M^1_{\generalring}$ to denote the set of all semi-finite signed measures on $\generalring$.

For each pair $A,B\in \generalring$ of disjoint sets, \ref{ZeroEmpty} and \ref{CountablyAdditiveList} combine to ensure that $\mu(A) + \mu(B)$ converges unconditionally to $\mu(A\cup B)$.
Therefore, an element of $\mathcal M^1_{\generalring}$ may take at most one of the values $\pm\infty$.  After all, suppose that $E,F\in \generalring$ are such that $\mu(E) = \infty$ and $\mu(F) = -\infty$.  Clearly
$\mu(E) + \mu(F\setminus E)$ converges unconditionally to $\mu(E\cup F)$, and hence, \ref{UC+} forces $\mu(E\cup F) = \infty$.  On the other hand, $\mu(F) + \mu(E\setminus F)$ converges unconditionally to $\mu(E\cup F)$,
implying that $\mu(E\cup F) = -\infty$, a contradiction.

Of course, there is extensive study of signed measures in the literature.  Importantly, if $\mu\in \mathcal M^1_{\generalring}$, then $\mu$ satisfies the required conditions of \cite[Definition 9.1.1]{Kesavan}, 
and is therefore a signed measure in the more traditional sense.  When we impose the stronger definition given above, we obtain a one-to-one correspondence between $\mathcal J^*$ and $\mathcal M^1_{\R}$.

\begin{thm} \label{MeasureIsomorphism}
	A map $\mu:\R\to [-\infty,\infty]$ is a charge on $\R$ if and only if it is a semi-finite signed measure on $\R$.  In particular, $\mathcal M_{\R}^1$ is a vector space over $\real$ and
	the map $\mu\mapsto c_\mu$ is a vector space isomorphism from $\mathcal M_{\R}^1$ to $\mathcal J^*$.
\end{thm}
\begin{proof}
	First assume that $\mu:\R\to \real$ is a charge.  Clearly $\mu$ satisfies properties \ref{ZeroEmpty} and \ref{CompactFinite}. Suppose that a set $A\in \R$ is expressed as a countable disjoint union
	\begin{equation*}
		A = \bigcup_{i=1}^\infty U_i,
	\end{equation*}
	where $U_i\in \R$ for all $i$.  Lemma \ref{RFacts}\ref{InfinitePairwiseDisjoint} implies that all but finitely many of the sets $U_i$ must be empty.  Therefore, the above union is really a finite union so that 
	\ref{CountablyAdditiveList} follows from \ref{FiniteAdditive}.  If $\mu:\R\to \real$ is a semi-finite signed measure, then clearly $\mu$ satisfies \ref{SimpleEmpty} and \ref{FiniteAdditive}.  Moreover, since every element of $\R$ is 
	compact (see Lemma \ref{RFacts}\ref{OpenCompact}), property \ref{CompactFinite} implies that $\mu$ takes real values.
	
	The second statement of the theorem follows from the fact that $\mathcal M_{\R}^0$ is a vector space over $\real$ along with Theorem \ref{SimpleMeasureIsomorphism}.
\end{proof}

The first statement of Theorem \ref{RingMeasure} is now a corollary to Theorem \ref{MeasureIsomorphism}, while the second statement follows from the fact that $f_\alpha$ is constant on the sets $Y(K,v)$
for any number field $K$ containing $\alpha$.  As a result, we have established the first of our two main results.

\section{Extensions of Consistent Maps to $\A$} \label{Algebra}

Let $\overline{\R} = \{B\subseteq Y: B^c\in \R\}$.  No element $A\in \overline{\R}$ can be compact as the collection $\{Y(\rat,p):p\in S\}$ is an open cover of $A$ having no finite subcover, and as a result, we know that 
$\R\cap \overline{\R} = \emptyset$.  Now define
\begin{equation*} 
	\A = \R\cup \overline{\R} = \left\{ A\subseteq Y: A\in \R\mbox{ or } A^c\in \R\right\}.
\end{equation*}
Theorem \ref{AlgebraMeasure} is an assertion about the smallest algebra of sets containing $\R$.  Our next lemma shows that $\A$ is indeed that algebra.

\begin{lem}\label{AAlgebra}
	$\A$ is the smallest algebra of sets on $Y$ containing $\R$.
\end{lem}
\begin{proof}
	We begin by showing that $\A$ is an algebra of sets.  Clearly $Y = Y\setminus \emptyset \in \overline{\R} \subseteq \A$.  To see that $\A$ is closed under unions, let $A,B\in \A$.
	If $A,B\in \R$ then since $\R$ is already known to be a ring of sets, we have $A\cup B\in \R\subseteq \A$.  Otherwise, assume without loss of generality that $B\in \overline{\R}$ and
	observe that $(A\cup B)^c = A^c \cap B^c$.  Since $A^c$ and $B^c$ are both closed and $B^c$ is compact, $(A\cup B)^c$ is a closed subset of a compact set, so it follows that $(A\cup B)^c$ is compact.
	Of course, $(A\cup B)^c$ is also open, so it follows from Lemma \ref{RFacts}\ref{OpenCompact} that $(A\cup B)^c\in \R$ as required.
		
	We now must prove that $\A$ is closed under set differences.  For this purpose, we again assume that $A,B\in \A$ and consider four cases.
	\begin{enumerate}[label={(\arabic*)}]
		\item If $A,B\in \R$ then since $\R$ is already known to be a ring of sets, we have $A\setminus B\in \R\subseteq \A$.
		\item If $A\in \R$ and $B\in \overline{\R}$ then we observe $A\setminus B = A\cap B^c$.  By our assumptions, $A$ and $B^c$ are both open and compact, and hence, $A\setminus B$
			is open and compact so it belongs to $\R$.
		\item If $A\in \overline{\R}$ and $B\in \R$ then we look at
			\begin{equation*}
				(A\setminus B)^c = (A\cap B^c)^c = A^c \cup B.
			\end{equation*}
			By our assumptions, both of the sets in this union belong to $\R$, and since $\R$ is a ring, we conclude that $(A\setminus B)^c\in \R$.  It now follows that $A\setminus B \in \overline{\R}$.
		\item If $A,B\in \overline{\R}$ then we once again observe that $A\setminus B = A\cap B^c$.  Both of these sets are clearly open, implying that $A\setminus B$ is open.  Additionally, $B^c$ is compact by Lemma 
			\ref{RFacts}\ref{OpenCompact}, so because $A$ and $B^c$ are closed, we obtain that $A\setminus B$ is compact.  Again applying Lemma \ref{RFacts}\ref{OpenCompact}, we conclude that $A\setminus B\in \R$.
	\end{enumerate}
	
	As every algebra must be closed under complementation, it clearly follows that $\A$ is the smallest algebra containing $\R$.
\end{proof}

Because of Theorem \ref{MeasureIsomorphism}, every consistent map $c:\mathcal J\to \real$ corresponds to a unique semi-finite signed measure $\mu\in \mathcal M^1_{\R}$ for which $\mu(Y(K,v)) = c(K,v)$.  
We recall from our introduction that $c:\R\to \real$ is used to denote this measure.  For each open subset $U$ of $Y$, a countable collection of non-empty pairwise disjoint sets $\Gamma\subseteq \R$ is called a {\it partition of $U$} if
\begin{equation*}
	U = \bigcup_{A\in \Gamma} A.
\end{equation*}
We emphasize that we are specifically requiring elements in a partition to belong to $\R$.  If we further have $\Gamma\subseteq Y(\mathcal J)$, then consistent with our earlier terminology, $\Gamma$ is called a {\it basis partition of $U$}.

Suppose that $\Gamma$ and $\Delta$ are partitions of the open subset $U$ of $Y$.  We say that $\Gamma$ is a {\it refinement} of $\Delta$ if, for each $B\in \Delta$, there exists a subset $\Gamma_B\subseteq \Gamma$ such that
\begin{equation*}
	B = \bigcup_{A\in \Gamma_B} A.
\end{equation*}
Since the sets in $\Gamma$ are assumed to be pairwise disjoint, each such set $A$ must be a subset of $B$ for some $B\in \Delta$.  We now establish a crucial property of refinements of partitions.

\begin{thm} \label{RefinementUnconditional}
	Suppose that $\Gamma$ and $\Delta$ are partitions of the same open subset $U$ of $Y$ and assume that $\Gamma$ is a refinement of $\Delta$.  Let $c:\mathcal J \to \real$ be a consistent map.
	If $\sum_{A\in \Gamma} c(A)$ converges unconditionally to a point $x \in[-\infty,\infty]$ then $\sum_{B\in \Delta} c(B)$ also converges unconditionally to $x$.
\end{thm}
\begin{proof}
	Since $\Gamma$ and $\Delta$ are both partitions of $U$, we clearly have.
	\begin{equation} \label{ZUnion}
		U = \bigcup_{A\in \Gamma} A = \bigcup_{B\in \Delta} B.
	\end{equation}
	Moreover, we have assumed that $\Gamma$ is a refinement of $\Delta$.  Hence, for each $B\in \Delta$, we may let $\Gamma_B\subseteq \Gamma$ be such that
	\begin{equation*}
		B = \bigcup_{A\in \Gamma_B} A.
	\end{equation*}
	Suppose now that $\sigma:\nat\to \Delta$ is a bijection.  To complete the proof, it is enough to show that $\sum_{n=1}^\infty c(\sigma(n))$ converges to $x$.
	
	To see this, we first write $\Gamma_n = \Gamma_{\sigma(n)}$ so that 
	\begin{equation} \label{SigmaUnion}
		\sigma(n) = \bigcup_{A\in \Gamma_n} A,
	\end{equation}
	and because $\Gamma$ is assumed to consist of pairwise disjoint sets, the right hand side of \eqref{SigmaUnion} is a pairwise disjoint union.  Moreover, every element of $\R$ is compact and no element of 
	$\Gamma$ may be empty, and hence, $\Gamma_n$ is a finite set.  As a result, we may write
	\begin{equation*}
		\Gamma_n = \{A(n,1),A(n,2),\ldots,A(n,m_n)\}
	\end{equation*}
	for some $m_n\in \nat$.  Then applying \eqref{SigmaUnion}, we obtain
	\begin{equation*}
		\sigma(n) = \bigcup_{i=1}^{m_n} A(n,i)\quad\mbox{and}\quad U = \bigcup_{n=1}^\infty \bigcup_{i=1}^{m_n} A(n,i),
	\end{equation*}
	where the second equality follows from \eqref{ZUnion} and the fact that $\sigma$ is a bijection.
	
	We now seek to define a bijection $\tau:\nat\to \Gamma$ so that we may apply our knowledge that $\sum_{A\in \Gamma} c(A)$ converges unconditionally to $x$.  To this end, we let
	\begin{equation*}
		M_n = \sum_{i=1}^n m_i
	\end{equation*}
	and note that $M_0 = 0$.  Then we define $\tau:\nat \to \Gamma$ by
	\begin{equation*}
		\tau(\ell) = A(n,\ell-M_{n-1}),
	\end{equation*}
	where $n$ is the unique non-negative integer such that $M_{n-1}< \ell \leq M_{n}$.  We claim that $\tau$ is a bijection.
	\begin{itemize}
		\item To establish surjectivity, we assume that $A\in \Gamma$.  Since $\Gamma$ is a refinement of $\Delta$ and $\sigma$ is surjective, there must exist $n\geq 1$ such that $A\subseteq \sigma(n)$.  This means that
			\begin{equation*}
				A\in \Gamma_{n} = \left\{ A(n,1),A(n,2),\ldots,A(n,m_{n})\right\}
			\end{equation*}
			and we may assume that $A = A(n,s)$, where $0< s < m_{n}$.  Letting $\ell = s +M_{n-1}$ we certainly have that $M_{n-1} < \ell \leq M_{n-1} + m_{n} = M_{n}$.  It now follows that
			$\tau(\ell) = A(n,\ell - M_{n-1}) = A(n,s)$, as required.
		\item To prove injectivity, we assume that $\tau(a) = \tau(b)$.  Further, we may let $M_{n-1} < a \leq M_{n}$ and $M_{n'-1} < b < M_{n'}$ which implies that $A(n,a-M_{n-1}) = A(n',b-M_{n'-1})$.
			If $n\ne n'$ then the injectivity of $\sigma$ means that $\sigma(n) \ne \sigma(n')$.  Since $\Delta$ is a collection of pairwise disjoint sets, we would obtain $\sigma(n)\cap \sigma(n') = \emptyset$
			contradicting the fact
			\begin{equation*}
				A(n,a-M_{n-1}) \subseteq \sigma(n)\quad\mbox{and}\quad A(n',b-M_{n'-1}) \subseteq \sigma(n').
			\end{equation*}
			We have now shown that $n = n'$ so that $A(n,a-M_{n-1}) = A(n,b-M_{n-1})$.  Then if $a\ne b$ then $A(n,a-M_{n-1})$ and $A(n,b-M_{n-1})$ would be distinct elements of $\Gamma_{n}$,
			also a contradiction.
		\end{itemize}
		
		We shall now prove that
		\begin{equation} \label{SeriesEquality}
			\sum_{\ell=1}^{M_N} c(\tau(\ell)) = \sum_{n=1}^N c(\sigma(n))\quad\mbox{for all } N\in \nat.
		\end{equation}
		To see this, observe that
		\begin{align*}
			\sum_{\ell=1}^{M_N} c(\tau(\ell)) & = \sum_{n=1}^N \sum_{\ell = M_{n-1} + 1}^{M_n} c(\tau(\ell)) \\
				& =  \sum_{n=1}^N \sum_{\ell = M_{n-1} + 1}^{M_n} c(A(n,\ell-M_{n-1})) \\
				& = \sum_{n=1}^N c\left( \bigcup _{\ell = M_{n-1} + 1}^{M_n} A(n,\ell-M_{n-1})\right),
		\end{align*}
		where the last equality follows from the fact that the sets $A(n,\ell-M_{n-1})$ are pairwise disjoint and $c$ is known to be a measure on $\R$.  Now simplifying this union, we conclude that
		\begin{align*}
			\sum_{\ell=1}^{M_N} c(\tau(\ell)) =  \sum_{n=1}^N c\left( \bigcup_{\ell = 1}^{m_n} A(n,\ell)\right)  = \sum_{n=1}^N c\left( \bigcup_{A\in \Gamma_n} A\right)  = \sum_{n=1}^N c\left(\sigma(n)\right),
		\end{align*}
		establishing \eqref{SeriesEquality}.
		
		As we have assumed $\sum_{A\in \Gamma} c(A)$ converges unconditionally to $x$, we know that $\sum_{n=1}^\infty c(\tau(\ell))$ converges to $x$.  In other words, the sequence
		\begin{equation*}
			\left\{ \sum_{n=1}^N c(\tau(\ell)) \right\}_{N=1}^\infty
		\end{equation*}
		converges to $x$.  Of course, the subsequence 
		\begin{equation*}
			\left\{ \sum_{n=1}^{M_N} c(\tau(\ell)) \right\}_{N=1}^\infty
		\end{equation*}
		must also converge to $x$ and the result follows from \eqref{SeriesEquality}.
\end{proof}

Since each element of $\mathcal R$ is a finite union of sets $Y(K,v)$, every partition of $U$ has a refinement that is a basis partition of $U$.  Moreover, if $\Gamma$ and $\Delta$ are arbitrary 
partitions of $U$ then the collection $\{A\cap B:A\in \Gamma,\ B\in \Delta\}$ is another partition of $U$ that is a refinement of both $\Gamma$ and $\Delta$.  These observations immediately yield 
the following observation as a consequence of Theorem \ref{RefinementUnconditional}.

\begin{cor} \label{ConvergenceEquivalences}
	If $U$ is an open subset of $Y$ and $c:\mathcal J\to \real$ is a consistent map then the following conditions are equivalent:
	\begin{enumerate}[label=(\roman*)]
		\item $\sum_{A\in \Gamma}c(A)$ converges unconditionally for all basis partitions $\Gamma$ of $U$.
		\item $\sum_{A\in \Gamma}c(A)$ converges unconditionally for all partitions $\Gamma$ of $U$.
		\item $\sum_{A\in \Gamma}c(A)$ converges unconditionally to the same point for all partitions $\Gamma$ of $U$.
	\end{enumerate}
\end{cor}

If the conditions of Proposition \ref{ConvergenceEquivalences} hold with $Y$ in place of $U$, then we recall that that $c$ is {\it globally consistent}.  We write $\overline{\mathcal J^*}$ to denote the set of globally consistent maps
and note that $\overline{\mathcal J^*}$ is a subspace of $\mathcal J^*$.  
Since each set $A\in \Gamma$ is compact, property \ref{CompactFinite} implies that $|c(A)|<\infty$ for all $A\in \Gamma$.  As a result, the definition of global consistency does not require us to consider \ref{UC+} and \ref{UC-}.
More importantly, global consistency is precisely the condition needed in order to extend a consistent map to a measure on $\A$.

\begin{thm} \label{ExtendedMeasureIsomorphism}
	The map $\mu\mapsto c_\mu$ is a bijection from $\mathcal M_{\A}^1$ to $\overline{\mathcal J^*}$.
\end{thm}
\begin{proof}
	As in the proof of Theorem \ref{MeasureIsomorphism}, we let $\phi$ be the map described in the theorem.  We must first prove that $\phi$ is well-defined.  Clearly $c_\mu$ is a consistent map, however, we must also prove that it
	is globally consistent.  To see this, we let $\Gamma$ be a basis partition of $Y$ so that $\mu(A) = c_\mu(A)$ for all $A\in \Gamma$.  However, since $\mu$ is a measure on $\A$ and $\cup_{A\in \Gamma} A \in \A$,
	countable additivity for $\mu$ implies that $\sum_{A\in \Gamma} \mu(A)$ converges unconditionally to a point in $[-\infty,\infty]$.  It follows that $\sum_{A\in \Gamma} c_{\mu}(A)$ converges unconditionally, as required.

	The injectivity of $\phi$ follows from Theorem \ref{MeasureIsomorphism}, so it remains only to prove surjectivity.  To this end, we assume that $c\in \overline{\mathcal J^*}$.   If $\Gamma$ is a partition of $Y$, 
	then the we define {\it index of $c$} by
	\begin{equation*}
		I(c) = \sum_{A\in \Gamma} c(A).
	\end{equation*} 
	Since $c$ is globally consistent, this definition does not depend on $\Gamma$.  Again by Theorem \ref{MeasureIsomorphism},
	there exists $\mu\in \mathcal M_{\R}^1$ such that $c_\mu = c$.  Now define $\nu:\A\to \real$ by
	\begin{equation*} \label{NuDef}
		\nu(A) = \begin{cases} \mu(A) & \mbox{if } A\in \R \\ 
						I(c) - \mu(A^c) & \mbox{if } A\in \overline{\R}.
				\end{cases}
	\end{equation*}
	We claim that $\nu$ is a measure on $\A$ such that $c_\nu = c$.
	
	If $A\in \overline{\R}$ then $A^c$ is a compact set in $\R$, and property \ref{CompactFinite} for $\mu$ implies that $\mu(A^c)$ is finite.  Therefore, the expression $I(c) - \mu(A^c)$ is well-defined
	even if $I(c) = \pm \infty$.  Now if $A\in Y(\mathcal J)$ then certainly $A\in \R$ so $c_{\nu}(A) = \nu(A) = \mu(A) = c_{\mu}(A) = c(A)$, implying that $c_\nu = c$.  Properties \ref{ZeroEmpty}
	and \ref{CompactFinite} follow trivially from the definition of $\nu$ so it remains only to establish countable additivity on $\A$.
	
	To this end, we assume that $\Delta$ is a countable collection of non-empty pairwise disjoint sets in $\A$.  If $A,B\in \overline{\R}$ then $(A\cap B)^c = A^c \cup B^c$ so that $(A\cap B)^c$ 
	is compact.  Therefore, there exists a finite set $Z$ of places of $\rat$ such that
	\begin{equation*}
		(A\cap B)^c \subseteq \bigcup_{p\in Z} Y(\rat,p).
	\end{equation*}
	There are clearly places of $\rat$ that are outside of $Z$, so $A\cap B$ is non-empty.  As a result, $\Delta$ contains at most one element of $\overline{\R}$.  For simplicity, we set
	\begin{equation*}
		F = \bigcup_{A\in \Delta} A
	\end{equation*}
	and assume that $F\in \A$.  We consider four cases:
	\begin{enumerate}[label={(\Roman*)}]
		\item Assume that $\Delta$ contains finitely many elements of $\R$ and no elements of $\overline{\R}$.  In this case, $F \in \R$ and we may apply the finite additivity
			property of $\mu$ on $\R$ to yield
			\begin{equation*}
				\nu\left(F\right) = \mu\left(F \right)  = \sum_{A\in \Delta} \mu(A) = \sum_{A\in \Delta} \nu(A).
			\end{equation*}
		\item Assume that $\Delta$ contains finitely many elements of $\R$ and one element $B\in \overline{\R}$.  This implies that $F \in \overline{\R}$.
			If $I(c) =\infty$ then $\nu(B) = \nu(F) = \infty$, so by \ref{UC+}, we know that $\sum_{A\in \Delta} \nu(A)$ converges unconditionally to $\infty$ and the required property holds.  
			A similar scenario occurs if $I(c) =-\infty$ so we may assume for the remainder of this proof that $I(c)\in \real$.
		
			Under this assumption, the definition of $\nu$ implies that $\nu(A)\in \real$ for all $A\in \A$ and $\sum_{A\in \Delta} \nu(A)$ is a finite sum of real numbers.
			Let $\Delta_0 = \Delta\setminus\{B\}$ and let $E = \bigcup_{A\in \Delta_0} A$ so we clearly have $E \in \R$ and $F = B\cup E$.   Then using the finite additivity of $\mu$ on $\R$, we obtain that
			\begin{align*}
				\sum_{A\in \Delta} \nu(A) & = \nu(B) + \sum_{A\in \Delta_0} \nu(A) \\
										& = I(c) - \mu(B^c) + \mu\left(E\right) \\
										& = I(c) - \left( \mu(B^c) - \mu\left(E\right) \right)
			\end{align*}
			The sets in $\Delta$ are assumed to be pairwise disjoint, and hence, $E\subseteq B^c$.  It now follows from finite additivity of $\mu$ as well as DeMorgan's Laws that
			\begin{align*} \label{StartingSum}
				\sum_{A\in \Delta} \nu(A) & = I(c) - \mu\left( B^c \setminus E\right) \\
										& = I(c) - \mu\left(B^c \cap E^c\right) \\
										& =I(c) - \mu\left( (B\cup E)^c\right) \\
										& = I(c) - \mu\left( F^c\right) \\
										& = \nu\left( F \right).
			\end{align*}
		\item Assume that $\Delta$ contains infinitely many elements of $\R$ and no elements of $\overline{\R}$.  In this situation, $F$ cannot be compact, so we have $F\in \overline{\R}$ and
			$\nu(F) = I(c) - \mu(F^c)$.  Letting $\Gamma = \Delta \cup\{F^c\}$, we clearly have that $\Gamma$ is a partition of $Y$ by elements of $\R$.  It follows that $\sum_{A\in \Gamma} c(A)$ converges 
			unconditionally to $I(c)$, and hence, $\sum_{A\in \Delta} c(A)$ converges unconditionally to $I(c) - c(F^c) = I(c) - \mu(F^c) = \nu(F)$.
		\item Assume that $\Delta$ contains infinitely many elements of $\R$ and one element $B\in \overline{\R}$.  We claim that these conditions lead to a contradiction.
			As before, let $\Delta_0 = \Delta\setminus\{B\}$ and let $E = \bigcup_{A\in \Delta_0} A$.  Clearly $F = B\cup E$ and $B\cap E = \emptyset$, so we obtain
			\begin{equation*}
				E = B^c \cap F.
			\end{equation*}
			Since $B\in \overline{\R}$, we know that $B^c\in \R$ and $B^c$ is open, closed and compact.  Additionally, $F\in \A$ so that $F$ is both open and closed.  It follows that
			$E$ is a closed subset of a compact set, and is therefore itself compact.  As a result, $E\in \R$ contradicting Lemma \ref{RFacts}\ref{InfinitePairwiseDisjoint}.
			
	\end{enumerate}
	
\end{proof}

Of course, Theorem \ref{AlgebraMeasure} now follows directly from Theorem \ref{ExtendedMeasureIsomorphism}.
It is worth noting that, even if $c$ is not globally consistent over $S$, then there are extensions of $\mu$ to a map $\nu:\mathcal A\to [-\infty,\infty]$ which satisfies \ref{ZeroEmpty} and \ref{CompactFinite}.
Indeed, one may select $r\in [-\infty,\infty]$ and define
\begin{equation*}
	\nu(A) = \begin{cases} \mu(A) & \mbox{if } A\in \R \\ 
						r - \mu(A^c) & \mbox{if } A\in \overline{\R}
				\end{cases}
\end{equation*}
to obtain such an extension, but $\nu$ will satisfy only finite additivity without necessarily having countable additivity.  In this case, countable additivity holds if and only if $\nu$ has continuity from below on $\mathcal A$.  
The proof of Theorem \ref{ExtendedMeasureIsomorphism} shows that continuity from below on $\mathcal A$ occurs precisely when $c$ is globally consistent and $r = I(c)$.

\bibliographystyle{abbrvnat}
\bibliography{MeasureCorrespondence}

\end{document}